\documentclass[11pt]{article}

\usepackage{amsmath}
\usepackage{amssymb}
\usepackage{authblk}
\usepackage{booktabs}
\usepackage{color}
\usepackage{epsfig,subfigure,stmaryrd,appendix}
\usepackage{epstopdf}
\usepackage{graphicx}
\usepackage{latexsym}
\usepackage{amsfonts}

\begin{document}  \baselineskip=16pt

\centerline{\LARGE How Many Numerical Eigenvalues can We Trust?}

\bigskip
\centerline{Zhimin Zhang \footnote{Beijing Computational Science Research Center, 100084, China;
and Department of Mathematics, Wayne State University, Detroit, MI 48202, USA.
This work is supported in part by
the US National Science Foundation under grant DMS-1115530.}}

\bigskip


{\it ABSTRACT. When using finite element and finite difference methods to approximate eigenvalues of $2m^{th}$-order elliptic problems,
the number of reliable numerical eigenvalues can be estimated in terms of the total degrees of freedom $N$ in resulting discrete systems.
The truth is worse than what we used to believe in that the percentage of reliable eigenvalues
decreases with an increased $N$, even though the number of reliable eigenvalues increases with $N$.}




\bigskip
\setcounter{equation}{0}
\setcounter{section}{1}

\noindent{\bf 1. Introduction}

When approximating PDE eigenvalue problems by numerical methods such as finite difference and finite element,
it is common knowledge that only a small portion of numerical eigenvalues are reliable. However, this knowledge is
only qualitative rather than quantitative in the literature \cite{babuska-osborn, strang-fix}. In this paper, we will
investigate the number of ``trusted" eigenvalues by the finite element
(and the related finite difference method results obtained from mass lumping)
approximation of $2m$th order elliptic PDE eigenvalue problems.
Our two model problems are the Laplace and bi-harmonic operators, for which a solid knowledge regarding
magnitudes of eigenvalues are available in the literature \cite{kac, levine-protter, li-yau, pleijel, polya, protter, weyl1}.
Combining this knowledge with {\it a priori} error estimates of the finite element method \cite{babuska-osborn, strang-fix},
we are able to figure out roughly how many ``reliable" eigenvalues can be obtained from numerical approximation under a
pre-determined convergence rate.

Let us begin with a simple example, which was used in \cite{babuska-osborn} and \cite{strang-fix}
for different purposes.
Approximating the one-dimensional eigenvalue problem
$$
-u'' = \lambda u, \quad u(0) = 0 = u(1); \qquad \lambda_j = (j\pi)^2, \quad u_j(x) = \sin(j\pi x).
$$
by linear finite element on the uniform mesh of $n$ subintervals results in an $(n-1)\times (n-1)$
linear algebraic system of generalized eigenvalue problems
\begin{equation}\label{1d}
\frac{1}{h^2}
\begin{pmatrix}
2 & -1    &      &  \\
-1&  2    &\ddots&  \\
  &\ddots&\ddots&-1\\
& & -1   &        2
\end{pmatrix} \vec{U} = \frac{\lambda}{6} \begin{pmatrix}
4 & 1    &      &  \\
1 & 4    &\ddots&  \\
  &\ddots&\ddots& 1\\
& &  1   &        4
\end{pmatrix} \vec{U},
\end{equation}
where $hn=1$.
Numerical eigenvalues are, for $\theta_j = j\pi h$, $j=1,2,\ldots,n-1$,
\begin{equation}\label{eig1d}
\lambda^h_j = \frac{3}{h^2} \frac{2(1-\cos\theta_j)}{2+\cos\theta_j}
= (j\pi)^2 + \frac{h^2}{12}(j\pi)^4 + \frac{h^4}{360}(j\pi)^6 + \cdots
\end{equation}
With the mass lumping (or equivalently the central finite difference scheme),
\begin{equation}\label{eig1dlump}
\lambda^{h,l}_j = \frac{2(1-\cos\theta_j)}{h^2} = (j\pi)^2 - \frac{h^2}{12}(j\pi)^4 + \frac{h^4}{360}(j\pi)^6 + \cdots
\end{equation}
Averaging the above two numerical eigenvalues $\lambda^{h*}_j = (\lambda^h_j+\lambda^{h,l}_j)/2$ yields
\begin{equation}\label{extrapl}
\lambda^{h*}_j = (j\pi)^2 + \frac{h^4}{360}(j\pi)^6 + \cdots
\end{equation}
We consider some special cases, the last and the middle (for even $n$) eigenvalues:
$$
\lambda^h_{n-1} = \frac{3}{h^2} \frac{2(1-\cos\theta_{n-1})}{2+\cos\theta_{n-1}} = 12n^2 - 9\pi^2 + O(\frac{1}{n^2}),
$$
$$
\lambda^{h,l}_{n-1} = \frac{2(1-\cos\theta_{n-1})}{h^2} = 4n^2 - \pi^2 + O(\frac{1}{n^2}),
$$
comparing with the exact value $\lambda_{n-1} = ((n-1)\pi)^2 = \pi^2n^2 - 2\pi n + \pi^2 \approx 9.8696 n^2 - 2\pi n + \pi^2$;
$$
\lambda^h_{n/2} = \frac{3}{h^2} = 3n^2, \quad \lambda^{h,l}_{n/2} = \frac{2}{h^2} = 2n^2,
$$
comparing with the exact value $\lambda_{n/2} = \displaystyle \frac{\pi^2}{4}n^2 \approx 2.4674 n^2$.
We see that in both cases, relative errors are of order $O(1)$.

Next, we investigate how many numerical eigenvalues have their relative errors converge in our expected rate.
Without loss of generality, let $s>1$ be a factor of $n$. From (\ref{eig1d}) and (\ref{eig1dlump}),
$$
\lambda^h_{n/s} = (\frac{n\pi}{s})^2 + \frac{h^2}{12}(\frac{n\pi}{s})^4 \pm \frac{h^4}{360}(\frac{n\pi}{s})^6 + \cdots,
$$
$$
\lambda^{h,l}_{n/s}  = (\frac{n\pi}{s})^2 + \frac{h^2}{12}(\frac{n\pi}{s})^4 - \frac{h^4}{360}(\frac{n\pi}{s})^6 + \cdots,
$$
Relative errors (note that $hn=1$):
\begin{equation}\label{relativ1}
\frac{\lambda^h_{n/s} - \lambda_{n/s}}{\lambda_{n/s}} = \frac{h^2}{12}(\frac{n\pi}{s})^2 + \frac{h^4}{360}(\frac{n\pi}{s})^4 \cdots
= \frac{\pi^2}{12}\frac{1}{s^2} + \frac{\pi^4}{360}\frac{1}{s^4} + \cdots,
\end{equation}
\begin{equation}\label{relativ2}
\frac{\lambda_{n/s} - \lambda^{h,l}_{n/s}}{\lambda_{n/s}}
= \frac{\pi^2}{12}\frac{1}{s^2} - \frac{\pi^4}{360}\frac{1}{s^4} + \cdots
\end{equation}

We observe in (\ref{relativ1}) and (\ref{relativ2}) that the convergence rate depends on $s$.
In order to have a quadratic convergence, as we usually expected from the standard convergence theory for eigenvalues,
we would need $s=O(n)$ and hence $n/s = O(1)$, which means that only the first few numerical eigenvalues qualify.
However, if we relax to linear convergence, we would require $s = O(\sqrt{n})$, which leads to $n/s = O(\sqrt{n})$,
i.e., we have about $\sqrt{n} = O(h^{-1/2})$ reliable numerical eigenvalues whose relative errors converge at least linearly.
If we only demand a weak convergence: $h^\alpha$ ($\alpha< 1$), we would have $s=O(n^{\alpha/2})$, $n/s = O(n^{1-\alpha/2})$,
hence more numerical eigenvalues would qualify.

Recall the extrapolated numerical eigenvalue (\ref{extrapl}), the relative error of $\lambda^{h*}_{n/s}$ is
$$
\frac{\lambda^{h*}_{n/s} - \lambda_{n/s}}{\lambda_{n/s}} = \frac{\pi^4}{360}\frac{1}{s^4} + O(s^{-6}).
$$
Following the same reasoning as in the above, for $\lambda^{h*}_j$,
we have roughly $\sqrt{n}$ numerical eigenvalues converge at least quadratically,
and about $n^{3/4}$ numerical eigenvalues
converge at least linearly, which is the same as under quadratic elements (which will be explained later by our general theory).

It is interesting to note that the relative error
$$
\frac{\lambda^{h*}_{n/2} - \lambda_{n/2}}{\lambda_{n/2}}
= \frac{2.5n^2 - n^2\pi^2/4}{n^2\pi^2/4}
= \frac{10}{\pi^2} - 1 \approx 1.32\%,
$$
is acceptable in practice.
However, for the last eigenvalue, the relative error for the averaged numerical approximation
$$
\frac{\lambda_{n-1} - \lambda^{h*}_{n-1}}{\lambda_{n-1}}
= \frac{(n-1)^2\pi^2 - 8n^2 + 5\pi^2 + O(n^{-2})}{(n-1)^2\pi^2}
\approx 1 - \frac{8}{\pi^2} \approx 18.94\%
$$
is not too much of an improvement over $\lambda^h_{n-1}$.

\bigskip


\setcounter{equation}{0}
\setcounter{section}{2}

\noindent{\bf 2. Model Problems and Finite Element Approximation}

\vskip.05in

Our model problems are the Laplace and bi-harmonic operators. The first model is
\begin{equation}\label{model1}
- \Delta u = \lambda u \quad\text{in}\quad \Omega, \qquad u = 0 \quad\text{on}\quad \partial\Omega.
\end{equation}
The countable sequence of eigenvalues is a classic result: $0 < \lambda_1 < \lambda_2 \le \cdots \le \lambda_k \le \cdots$
tending to $+\infty$ and a sequence of corresponding eigenfunctions $u_1, u_2, \ldots, u_k, \ldots$
such that each $u_k$ satisfies (\ref{model1}).
The eigenfunctions are orthogonal in $L_2(\Omega)$ and they are customarily normalized so
that $\|u_k\|_{L_2(\Omega)} = 1$ for all $k$. The second model is
\begin{equation}\label{model2}
\Delta^2 u = \mu u \quad\text{in}\quad \Omega, \qquad u = 0, \;\; \frac{\partial u}{\partial \pmb n} = 0 \quad\text{on}\quad \partial\Omega,
\end{equation}
 which also has a set of normal modes of vibration $0 < \mu_1 \le \mu_2 \cdots \le \mu_k \le \cdots$ and a corresponding set of eigenfunctions.

\vskip.1in

Next, we consider finite element approximations of eigenvalue problems of elliptic operators.
First, we introduce some notations and classical results starting from the following abstract form:
Find $(\lambda, u) \in R\times H$, for a suitable Hilbert space $H$, such that
\begin{equation}\label{eig-modle}
a(u,v) = \lambda b(u,v), \quad\forall v \in H,
\end{equation}
where $a(\cdot,\cdot)$ is a bilinear form on $H$ deduced from a self-adjoint $2m$-order elliptic operator on a smooth domain in $R^d$ and
$b(\cdot,\cdot)$ is a symmetric positive definite bilinear form on $H$.

In our discussion, we include multiple eigenvalues.
We consider conforming finite element methods.
 To be more precise, let $(\lambda^h,u_h)$ solve (\ref{eig-modle}) on a finite dimensional subspace $V_h\in H$,
 which consists of polynomials of degrees fewer than
  $k$ with $k> m\ge 1$ on a shape regular and quasi-uniform triangulation ${\cal T}_h$.
 We further define the energy norm $\|\cdot\|_H$ on $H$ and introduce the following notations \cite{babuska-osborn}:
$$
M(\lambda) = \{ u: \; \text{all eigenvectors of (\ref{eig-modle}) corresponding to}\; \lambda, \|u\|_H = 1\},
$$
$$
\epsilon_h(\lambda) = \sup_{u\in M(\lambda)} \inf_{\chi\in V_h} \|u-\chi\|_H.
$$

\bigskip


\setcounter{equation}{0}
\setcounter{section}{3}

\noindent{\bf 3. Theoretical Issues}

\vskip.05in

More than a century ago,
Weyl obtained the following asymptotic behavior for the $n$th eigenvalue of the Laplace operator, problem (\ref{model1}) \cite{weyl1}:
\begin{equation}\label{weyl}
\lambda_n \approx 4\pi^2 \left( \frac{n}{\omega_d|\Omega|} \right)^{2/d},
\end{equation}
where $\omega_d = \pi^{d/2}/\Gamma(1+d/2)$ denote the volume of the unit ball in $R^d$. Pleijel extended it to (\ref{model2}) by showing that \cite{pleijel}
\begin{equation}\label{plei}
\mu_n \approx 16\pi^4 \left( \frac{n}{\omega_d|\Omega|} \right)^{4/d}.
\end{equation}
Further, Polya proved in 1961 that for a tiling domain $\Omega$ \cite{polya}
\begin{equation}\label{pol}
 \lambda_n \ge 4\pi^2 \left( \frac{n}{\omega_d|\Omega|} \right)^{2/d},
\end{equation}
and he further conjectured that (\ref{pol}) is valid for arbitrary domains. So far, the best known result along this line was due to Li-Yau \cite{li-yau}
\begin{equation}\label{ly}
\sum_{i=1}^n \lambda_i \ge \frac{dn}{d+2} 4\pi^2 \left( \frac{n}{\omega_d|\Omega|} \right)^{2/d} \quad
\to \lambda_n \ge \frac{d}{d+2} 4\pi^2 \left( \frac{n}{\omega_d|\Omega|} \right)^{2/d}.
\end{equation}
For more details, the reader is referred to \cite{kac, protter} and references therein. For the purpose of this article, (\ref{weyl}) and (\ref{plei}) would be sufficient, and it is well know that $\mu_k \ge \lambda_k^2$.


{\bf Proposition 1} (Strang-Fix 1973 \cite{strang-fix}): Using conforming finite element spaces of degree smaller than
$k (> m \ge 1)$ on a shape regular and quasi-uniform triangulation
in approximating a $2m$-order elliptic operator eigenvalue problem, the following relative error bound holds for the $i$th eigenvalue,
\begin{equation}\label{stran}
0 < \frac{\lambda^h_i - \lambda_i}{\lambda_i} \le C h^{2(k-m)}\lambda^{k/m-1}_i,
\end{equation}
when associated eigenfunction is sufficiently smooth.

\vskip.1in

{\bf Proposition 2} (Babu\v{s}ka-Osborn 1989 \cite{babuska-osborn}): Under the same assumption as in Proposition 1, for the $i$th eigenfunction,
there exist positive constants $c$ and $C$ independent of $(\lambda_i, u_i)$ and $h$, such that
\begin{equation}\label{babuska}
c\epsilon_h^2(\lambda) \le \lambda^h_i-\lambda_i \le C\epsilon_h^2(\lambda).
\end{equation}

{\bf Proposition 3} \cite{lin-xie-xu}: Under the same assumption as in Proposition 1, for the $i$th eigenfunction,
there exist positive constants $C$ independent of $(\lambda_i, u_i)$ and $h$, such that
$$
\epsilon_h(\lambda_i) \ge Ch^{k-m}\lambda^{k/(2m)}_i.
$$

Combining Propositions 1-3, we conclude that the relative error for $i$th eigenvalue behaves like
\begin{equation}\label{eig-err}
\frac{\lambda^h_i - \lambda_i}{\lambda_i} \approx h^{2(k-m)}\lambda^{k/m-1}_i.
\end{equation}
Here ``$\approx$" means both upper and lower error bounds are of the same order.
This relative error estimate is the first basic assumption for our following main theorem.
Our second assumption is about the asymptotic growth of the exact eigenvalues
\begin{equation}\label{eig-growth}
\lambda_j = O(j^{2m/d}).
\end{equation}


{\bf Theorem}. Suppose that we solve a $2m$-order elliptic equation on a domain $\Omega \subset R^d$ by the finite element method
(conforming or non-conforming) of polynomial degree $k-1$ under a shape regular and quasi-uniform mesh with mesh-parameter $h$.
Assume that the exact eigenvalues grow as (\ref{eig-growth}) and the relative error can be estimated by (\ref{eig-err}).
Then there are about
\begin{equation}\label{jn}
j_N = N^{(k-m-\alpha/2)/(k-m)}(k-1)^{-d(k-m-\alpha/2)/(k-m)}
\end{equation}
reliable numerical eigenvalues with relative error of $\lambda_{J_N}$, converges at rate $h^\alpha$ for $\alpha \in (0, 2(k-m)]$.
Here $N$ is the total degrees of freedom.

Proof:
If we want the relative error (\ref{eig-err}) to converge at rate $h^\alpha$, we would have
\begin{equation}\label{id1}
h^{2(k-m)}\lambda^{k/m-1}_j = h^\alpha, \quad\text{or}\quad  (h\lambda_j^{1/(2m)})^{2(k-m)} = h^\alpha.
\end{equation}
According to (\ref{eig-growth}), identity (\ref{id1}) leads to
\begin{equation}\label{id2}
(hj^{1/d})^{2(k-m)} = h^\alpha, \quad\text{or}\quad j = h^{-d(k-m-\alpha/2)/(k-m)}.
\end{equation}
On a shape regular and quasi-uniform mesh with mesh parameter $h$ in $R^d$, a piecewise polynomial space of degree $k-1$
has the total degrees of freedom $N = O((k-1)^dh^{-d})$, therefore there are about
\begin{equation}\label{id1a}
j_N = N^{(k-m-\alpha/2)/(k-m)}(k-1)^{-d(k-m-\alpha/2)/(k-m)}
\end{equation}
 reliable numerical eigenvalues with $\lambda_{J_N} $ converges at rate $h^\alpha$.

\vskip.1in

{\it Remark 1}.
From (\ref{weyl}) and (\ref{plei}), we know (\ref{eig-growth}) is valid at least for
the Laplace ($m=1$) and bi-harmonic ($m=2$) operators. From Propositions 1-3,
our theory covers standard conforming finite elements.
In addition, our theory would cover nonconforming finite elements, even discontinuous Galerkin methods, as long as (\ref{eig-err}) is valid.


\vskip.05in

1) Special case, second-order problem $m=1$.

1.1) Linear element $k=2$ "
(For notational consistency with \cite{strang-fix}, $k=2$ represent linear element here.),
there are about $N^{1-\alpha/2}$
reliable numerical eigenvalues converging at least with rate $h^\alpha$.
If we demand optimal convergence rate $\alpha = 2$, we have $N^0=1$,
which indicates that only some earlier eigenvalues can be approximated at quadratic rate $h^2$.
However, if we relax the convergence requirement to linear rate $\alpha=1$,
we would have about $\sqrt{N}$ numerical eigenvalues qualify.

1.2) For quadratic element $k=3$,
there are about $2^{-d(1-\alpha/4)}N^{1-\alpha/4}$ reliable numerical eigenvalues converging at least with rate $h^\alpha$.
If we demand optimal convergence rate $\alpha = 4$,
we have $N^0=1$, which indicates that only some earlier eigenvalues can be approximated at quartic rate $h^4$.
However, if we demand only second order convergence with $\alpha = 2$, we would have about $2^{-d/2}\sqrt{N}$ numerical eigenvalues qualify.
If we further relax the convergence rate to linear, there will be roughly $2^{-3d/4}N^{3/4}$ qualified numerical eigenvalues.

\vskip.05in

2) Special case, fourth-order problem $m=2$. Identity (\ref{id1a}) leads to
$$j_N = N^{(k-2-\alpha/2)/(k-2)}(k-1)^{-d(k-2-\alpha/2)/(k-2)}.$$
In order to have a meaningful rate of convergence, $k>2$, which means that we need at least quadratic element.
Indeed, a fourth-order problem requires higher order elements.

2.1) Quadratic element $k=3$, there are about $N^{1-\alpha/2}2^{-d(1-\alpha/2)}$ reliable numerical eigenvalues converging at least with rate $h^\alpha$.
If we demand optimal convergence rate $\alpha = 2$, we have $N^0=1$,
which indicates that only some earlier eigenvalues can be approximated at quadratic rate $h^2$.
However, if we relax the convergence requirement to linear rate $\alpha=1$,
we would have about $2^{-d/2}\sqrt{N}$ numerical eigenvalues qualify.

2.2) Cubic element $k=4$, there are about $3^{-d(1-\alpha/4)}N^{1-\alpha/4}$ reliable numerical eigenvalues converging at least with rate $h^\alpha$.
If we demand optimal convergence rate $\alpha = 4$, we have $N^0=1$,
which indicates that only some earlier eigenvalues can be approximated at quatic rate $h^4$.
However, if we demand only second order convergence with $\alpha = 2$,
we would have about $3^{-d/2}\sqrt{N}$ numerical eigenvalues qualify.
If we further relax the convergence rate to linear, there will be $3^{-3d/4}N^{3/4}$ qualified numerical eigenvalues.





\bigskip


\setcounter{equation}{0}
\setcounter{section}{4}

\noindent{\bf 4. Further Examples and Discussions}

\vskip.05in

In this section, we present some further examples to illustrate our general theory.
We apply both linear and bi-linear elements ($k=2$) to
 the Laplacian ($m=1$) eigenvalue problem (\ref{model1}) when $\Omega$ is the unit square ($d=2$).

\noindent 2.1. Linear element on regular triangulation with $(n-1)^2$ interior nodes.
The resulting linear algebraic generalized eigenvalue problem is:
\begin{equation}\label{linear}
\frac{1}{h^2}
\begin{pmatrix} A & -I &  &  \\ -I & A & -I &  \\ & \ddots & \ddots & \ddots \\ &  & -I & A \end{pmatrix}
\vec{U} = \frac{\lambda}{12}
\begin{pmatrix} H  & D &  &  \\ D' & H & D &  \\ & \ddots & \ddots & \ddots \\ &  & D' & H \end{pmatrix} \vec{U},
\end{equation}
where $I$ is the $(n-1)\times (n-1)$ identity matrix.
$$
A = \begin{pmatrix} 4 & -1 & & \\ -1 & 4 & \ddots & \\ & \ddots & \ddots & -1 \\ & & -1 & 4 \end{pmatrix}, \quad
H = \begin{pmatrix} 6 & 1 & & \\ 1 & 6 & \ddots & \\ & \ddots & \ddots & 1 \\ & & 1 & 6 \end{pmatrix}, \quad
D = \begin{pmatrix} 1 &  & & \\  1 & 1 &  & \\ & \ddots & \ddots &  \\ & &  1 & 1 \end{pmatrix}.
$$
By mass lumping, the mass matrix becomes an identity matrix,
which is equivalent to the counterpart 5-point finite difference scheme.
The eigenvalues of $A$ are $4-2\cos\theta_j$.
Therefore, numerical eigenvalues with mass lumping are
\begin{equation}\label{eiglinear}
\lambda^{h,l}_{j,k} = \frac{1}{h^2} ( 4-2\cos\theta_j - 2\cos\theta_k ).
\end{equation}
Recall the one dimensional linear element with mass lumping case (\ref{eig1dlump}), we see that
$$
\lambda^{h,l}_{j,k} = \lambda^{h,l}_j + \lambda^{h,l}_k.
$$
Therefore, we have the same conclusion as under one dimension for the two extremal cases
$\lambda^{h,l}_{n-1,n-1}$ and $\lambda^{h,l}_{n/2,n/2}$; in addition, we have, for $k, j << n$,
\begin{equation}\label{eigl}
\lambda^{h,l}_{j,k} = (j\pi)^2 + (k\pi)^2 - \frac{h^2}{12} [(j\pi)^4 + (k\pi)^4] + \frac{h^4}{360} [(j\pi)^6 + (k\pi)^6] + \cdots
\end{equation}
Here the total degree of freedom $N = (n-1)^2$, in which case $\sqrt{N} \approx n$.

In the two dimensional setting, we order eigenvalues according to their magnitudes.
There are three natural ways to do that.

1) Triangular ordering: for each level $\ell=1,2,\ldots$,
we choose $\lambda_{1,\ell}, \lambda_{\ell,1}, \lambda_{2,\ell-1}, \lambda_{\ell-1,2}, \ldots$,
the last one is $\lambda_{\ell/2,\ell/2}$ when $\ell$ is even, and the last two are
$\lambda_{(\ell-1)/2,(\ell+1)/2}, \lambda_{(\ell+1)/2,(\ell-1)/2}$ when $\ell$ is odd.
There are $(\ell+1)(\ell+2)/2$ eigenvalues in the triangle, and the two with the largest magnitude
are $\lambda_{1,\ell}$ and $\lambda_{\ell,1}$.

2) Square ordering: for each level, we choose
$\lambda_{1,\ell}, \lambda_{2,\ell}, \ldots, \lambda_{\ell,\ell}, \lambda_{\ell,\ell-1}, \ldots, \lambda_{\ell,2}, \lambda_{\ell,1}$.
The one with the largest magnitude is $\lambda_{\ell,\ell}$ and there are $\ell^2$ eigenvalues in the square.

3) Circular ordering: draw quarter-circles in the first quadrant with radiuses $1+1, 1+2^2, 1+3^2, \ldots, 1+\ell^2$,
group lattice points in each ring formed by adjacent circles. There are about $\pi\ell^2/2$ eigenvalues in the quarter disk.
On the outer quarter circle, in addition to $\lambda_{1,\ell}$ and $\lambda_{\ell,1}$ (with radius $1+\ell^2$),
there might be several other eigenvalues that have the largest magnitude, depending on how many lattice points on the circle.

Note that in the above three ordering strategies, the number of eigenvalues differ by factors less than 2.
For simplicity, we use the second strategy, square ordering, in which case,
the first $O(\sqrt{N})$ numerical eigenvalues are then obtained by setting $\ell = \sqrt{n}$,
and their relative errors are less than
$$
\frac{\lambda_{\sqrt{n},\sqrt{n}} - \lambda^{h,l}_{\sqrt{n},\sqrt{n}}}{\lambda_{\sqrt{n},\sqrt{n}}}
=  \frac{h^2}{12} \frac{(\sqrt{n}\pi)^4 + (\sqrt{n}\pi)^4}{(\sqrt{n}\pi)^2 + (\sqrt{n}\pi)^2} + O(h^2)
= \frac{h\pi^2}{12} + O(h^2).
$$
Here, we assume that $n$ is a complete square or otherwise we can round it to its nearest complete square without loss of its order.

\vskip.05in

{\it Remark 2}. We see that there are about $O(\sqrt{N}) = O(n) = O(h^{-1})$ numerical eigenvalues with relative errors converging linearly.
This conclusion is the same as the one dimensional case if measured by the total degrees of freedom $N$. Indeed, it is
also true for any dimension if linear element is used.
However, if measured by the mesh size $h$, the conclusion would be different
for linear element, in which case there are about $O(h^{-1/2})$, $O(h^{-1})$, and $O(h^{-3/2})$ numerical eigenvalues, which have relative errors converging linearly for the one, two, and three dimensional situations.

\vskip.1in

\noindent 2.2. Bilinear element on $n\times n$ square partition. The resulting linear algebraic system is:
\begin{equation}\label{bilinear}
\frac{1}{3h^2} \begin{pmatrix}
A & -B    &      &  \\
-B & A    &\ddots&  \\
  &\ddots&\ddots& -B\\
& &  -B   &        A
\end{pmatrix} \vec{U} =
\frac{\lambda}{36} \begin{pmatrix}
4C &  C   &      &   \\
 C & 4C   &\ddots&   \\
   &\ddots&\ddots&  C\\
   &      & C    & 4C
\end{pmatrix} \vec{U},
\end{equation}
where
$$
A = \begin{pmatrix}
8 & -1    &      &  \\
-1&  8    &\ddots&  \\
  &\ddots&\ddots&-1\\
& & -1   &        8
\end{pmatrix}, \quad
B = \begin{pmatrix}
1 & 1    &      &  \\
1 & 1    &\ddots&  \\
  &\ddots&\ddots& 1\\
& &  1   &        1
\end{pmatrix}; \quad
C = \begin{pmatrix}
4 & 1    &      &  \\
1 & 4    &\ddots&  \\
  &\ddots&\ddots& 1\\
& &  1   &        4
\end{pmatrix}.
$$

We first consider mass lumping, which results in the identity mass matrix and therefore is equivalent to the 9-point finite difference scheme.
Note that eigenvalues of $A$ and $B$ are $8-2\cos\theta_j$ and $1+2\cos\theta_k$, respectively.
Then numerical eigenvalues for the mass lumping scheme are
\begin{eqnarray}\label{eigbilump}
\lambda^{h,lb}_{j,k} &=& \frac{1}{3h^2} [8-2\cos\theta_j - 2(1+2\cos\theta_j)\cos\theta_k] \nonumber\\
&=& \frac{1}{3h^2} [ 2(1-\cos\theta_j) + 2(1-\cos\theta_k) + 4(1-\cos\theta_j\cos\theta_k) ].
\end{eqnarray}
When $k,j<<n$, using
$$
4(1-\cos\theta_j\cos\theta_k) = 2(\theta_j^2+\theta_k^2) - \frac{1}{6}(\theta_j^4+\theta_k^4) - \theta_j^2\theta_k^2 + \cdots
$$
we obtain
\begin{eqnarray}\label{eigblump}
\lambda^{h,lb}_{j,k} &=& (j\pi)^2 + (k\pi)^2 - \frac{h^2}{12} [(j\pi)^4 + (k\pi)^4 + 4(j\pi)^2(k\pi)^2] \nonumber\\
&& \; + \frac{h^4}{120} [ (j\pi)^6 + (k\pi)^6 + 10(j\pi)^4(k\pi)^2 + 10(j\pi)^2(k\pi)^4 ] + \cdots
\end{eqnarray}
Next, we consider the original bilinear case. Note that eigenvalues of the mass matrix on the right hand side of (\ref{bilinear}) are
$$
\frac{1}{36} [ (4+2\cos\theta_j) + 2(4+2\cos\theta_j)\cos\theta_k ]
= \frac{1}{9} (4 + 2(\cos\theta_j+\cos\theta_k) + \cos\theta_j\cos\theta_k).
$$
Therefore, eigenvalues from the bilinear element are
\begin{eqnarray}\label{eigbi}
\lambda^{h,b}_{j,k} &=& \frac{3}{h^2} \frac{2(1-\cos\theta_j) + 2(1-\cos\theta_k) + 4(1-\cos\theta_j\cos\theta_k) }
{4 + 2(\cos\theta_j+\cos\theta_k) + \cos\theta_j\cos\theta_k} \nonumber\\
&=& \frac{1}{h^2} ( \theta^2_j + \theta^2_k ) + \frac{1}{12} (\theta^4_j + \theta^4_k) + \frac{1}{360} (\theta^6_j + \theta^6_k) - \frac {17}{60480}(\theta^8_j + \theta^8_k) + \cdots \nonumber\\
&=& (j\pi)^2 + (k\pi)^2 + \frac{h^2}{12} ( (j\pi)^4 + (k\pi)^4 ) + \frac{h^4}{360} ( (j\pi)^6 + (k\pi)^6 ) \nonumber\\
&& \; - \frac {17h^6}{60480}((j\pi)^8 + (k\pi)^8) + \cdots.
\end{eqnarray}
Let us consider again some special cases. From (\ref{eigbilump}), we have
\begin{eqnarray*}
\lambda^{h,lb}_{n-1,n-1} &=& \frac{1}{3h^2} [ 4(1+\cos\frac{\pi}{n}) + 4(1-\cos^2\frac{\pi}{n}) ] \\
&=& \frac{4}{3h^2} (1+\cos\frac{\pi}{n})(2-\cos\frac{\pi}{n}) \\
&=& \frac{4}{3h^2} (2-\frac{1}{2}(\frac{\pi}{n})^2 + \cdots)(1+\frac{1}{2}(\frac{\pi}{n})^2 + \cdots) \\
&=& \frac{4}{3h^2} (2+\frac{1}{2}(\frac{\pi}{n})^2 + \cdots) = \frac{8}{3}n^2 + \frac{2\pi^2}{3} + O(n^{-2});
\end{eqnarray*}
and from (\ref{eigbi}), we obtain
\begin{eqnarray*}
\lambda^{h,b}_{n-1,n-1} &=& \frac{3}{h^2} \frac{ 4(1+\cos\frac{\pi}{n})(2-\cos\frac{\pi}{n}) } { 4-4\cos\frac{\pi}{n}+\cos^2\frac{\pi}{n} }
= \frac{12}{h^2} \frac{ 1+\cos\frac{\pi}{n} } { 2-\cos\frac{\pi}{n} } \\
&=& \frac{12}{h^2} (2-\frac{1}{2}(\frac{\pi}{n})^2 + \cdots)(1+\frac{1}{2}(\frac{\pi}{n})^2 + \cdots)^{-1} \\
&=& \frac{12}{h^2} (2-\frac{3}{2}(\frac{\pi}{n})^2 + \cdots) = 24n^2 - 18\pi^2 + O(n^{-2}).
\end{eqnarray*}
Comparing with the exact eigenvalue $\lambda_{n-1,n-1} = 2(n-1)^2\pi^2$, we see that the relative errors for bilinear element with and without lumping are of order $O(1)$.
$$
\lambda^{h,lb}_{n/2,n/2} = \frac{8}{3}n^2,  \quad \lambda^{h,b}_{n/2,n/2} = 6n^2. \qquad (\lambda_{n/2,n/2} = \frac{\pi^2}{2} n^2).
$$
Again, their relative errors are $O(1)$.

\vskip.05in

Recall the one dimensional case, by extrapolating $\lambda^h_j$ and $\lambda^{h,l}_j$, we obtain about $O(N^{3/4})$ numerical eigenvalues whose relative errors converge linearly. In the two dimensional case, this extrapolation is done between $\lambda^{h,l}_{j,k}$ (\ref{eigl}) and $\lambda^{h,b}_{j,k}$ (\ref{eigbi}). Indeed,
\begin{equation}\label{extrap2}
\lambda^{h*}_{j,k} = \frac{1}{2} ( \lambda^{h,b}_{j,k} + \lambda^{h,l}_{j,k} )
= (j\pi)^2 + (k\pi)^2 + \frac{h^4}{360}[ (j\pi)^6 + (k\pi)^6 ] + \cdots
\end{equation}
Comparing with (\ref{extrapl}), what a match! It is straightforward to calculate
$$
\lambda^{h,*}_{n-1,n-1} = 16n^2 - 10\pi^2 + O(n^{-2}), \quad \lambda^{h*}_{n/2,n/2} = 5n^2.
$$
Therefore,
$$
\frac{\lambda_{n-1,n-1} - \lambda^{h,*}_{n-1,n-1}}{\lambda_{n-1,n-1}}
= \frac{2(n-1)^2\pi^2 - (16n^2 - 10\pi^2 + O(n^{-2}))}{2(n-1)^2\pi^2}
 \approx 1- \frac{8}{\pi^2},
$$
$$
\frac{\lambda^{h,*}_{n/2,n/2} - \lambda_{n/2,n/2}}{\lambda_{n/2,n/2}}
= \frac{5n^2 - n^2\pi^2/2}{n^2\pi^2/2}
= \frac{10}{\pi^2} - 1.
$$
which are the same as in the one dimensional situation.

{\it Remark 3}. We see that in the two dimensional situation, linear element with mass lumping plays the same role as in the one dimensional case;
However, it is the bilinear element (not the linear element) that plays relatively the same role as linear element for the one-dimensional case.

\bigskip


\setcounter{equation}{0}
\setcounter{section}{5}

\noindent{\bf 5. A Comparison with the Spectral Method}

\vskip.1in

In this section, we compare finite element methods with spectral methods on the most simple setting in both one and two dimensional situations.
According to Weideman-Trefethen \cite{weideman-trefethen}, there are about $2/\pi$ portion ``trusted" eigenvalues for the polynomial spectral method
in 1-D. We understand ``trusted" to mean at least $O(N^{-1})$ accuracy with polynomial degree $N$.

{\bf Example 1}. Consider eigenvalue problem $-u'' = \lambda u$ on $[-1,1]$ with $u(-1) = 0 = u(1)$.
We divide the interval $[-1,1]$ into $2^{13} = 8192$ equal length subintervals and apply linear finite element method. With the mass lumping,
we obtain the central finite difference scheme.
We then apply Legendre spectral method with polynomial degree $N=8192$. Note that for the Legendre spectral method,
the stiffness matrix is diagonal and the mass matrix is 5-diagonal with only $3N-2$ none-zero entries, and condition number is not an issue.
In Figure 1, we depict relative errors for 8191 eigenvalues for all three methods and draw a horizontal line $y=0.012\% \approx 1/8192$ and two vertical lines
$x=\sqrt{N}$ and $x=2N/\pi$. We see that the error curves of finite element, finite difference methods pass through the intersection of
$y=0.012\%$ and $x=\sqrt{N}$, and the error curve of the spectral method cuts the intersection of $y=0.012\%$ and $2N/\pi$.
This observation is inconsistent with our theoretical prediction that there are about $O(\sqrt{N})$ (for $N=8192$, this is about $1\%$)
numerical eigenvalues that have their relative errors converge linearly.
Note that for $h=N^{-1}$, linear convergence rate means that the relative error is about $0.012\%$.
Our calculation also confirms that there are about $2/\pi$
numerical eigenvalues for which the relative error converges at rate $O(N^{-1})$ for Legendre polynomial spectral method.

\vskip.1in

{\bf Example 2}. Consider the eigenvalue problem (\ref{model1}) when $\Omega = [-1,1]^2$.
We use four different numerical methods to solve it:
1) Divide $\Omega$ into $2^{12} \times 2^{12}$ equal sub-squares and apply linear finite element method;
2) Introduce mass lumping to the linear finite element method and obtain the 9-point finite difference method;
3) use quadratic element on the $2^{11} \times 2^{11}$ equal sub-squares partition;
and 4) apply Legendre spectral method with polynomial degree $2^{12}=4096$ in each direction.

All four methods have the same total degrees of freedom $N = 2^{24}$.
Note that for the Legendre spectral method,
the stiffness matrix is again diagonal and the mass matrix is block 5-diagonal with only $(3N-2)^2$ none-zero entries.
The condition number in 2D is also not an issue.

In Figure 2, we draw $0.024\% \approx 1/4096$ relative error regions for each of the four methods.
Inside the region, the relative error is less than $0.024\%$.
We see that only $0.03\% (\approx 1.9N^{-1/2})$ numerical eigenvalues from linear finite element
or the 9-point finite difference method qualify,
and about $1.15\% (\approx 9.2N^{-1/4})$ numerical eigenvalues from quadratic finite element qualify.
The percentage of values that qualify increases to $40.58\% (\approx (2/\pi)^2)$ for the Legendre spectral method.

\begin{figure}
\begin{center}
\includegraphics{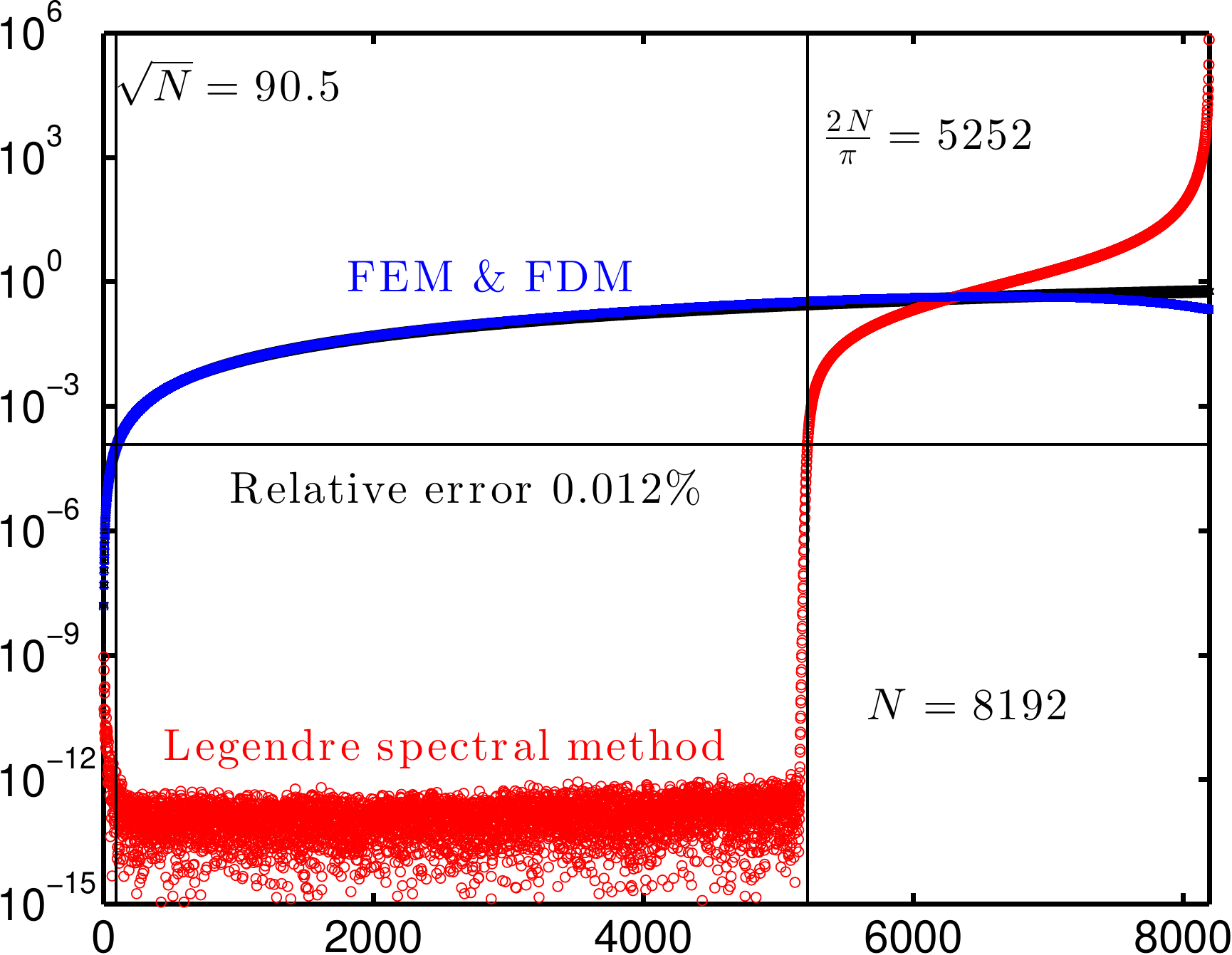}
\end{center}
\end{figure}

\begin{figure}
\begin{center}
\includegraphics{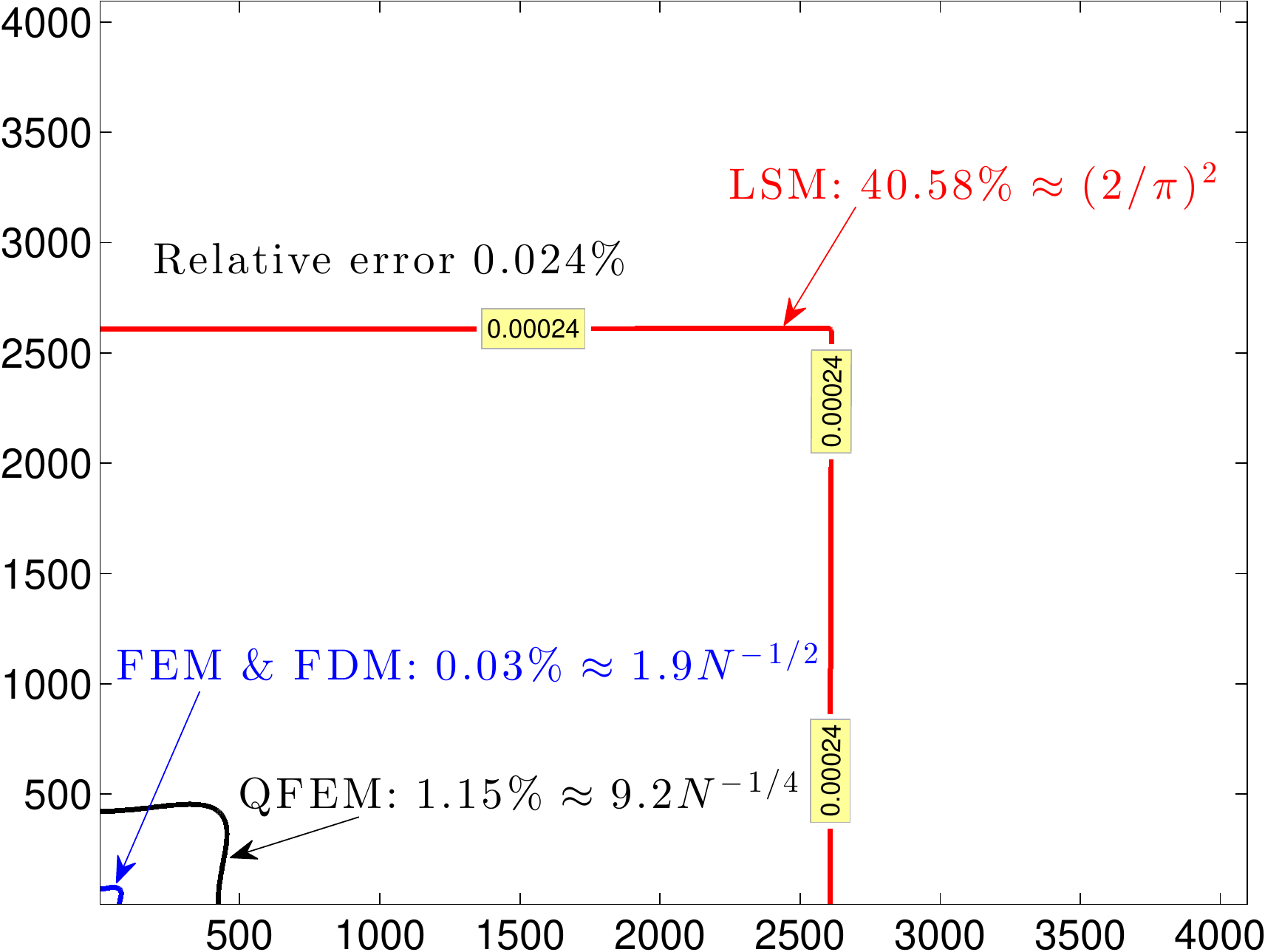}
\end{center}
\end{figure}

\bigskip

{\bf Conclusion Remarks}. From the above discussion, we see that only a small portion of numerical eigenvalues 
obtained from finite element (finite difference) methods are reliable even under the most favorable situation,
i.e., eigenfunctions are sufficiently smooth and round-off errors are taking off the picture. 
Although the number of reliable eigenvalues increases with an increased computational scale $N$,
the percentage of reliable eigenvalues (compared with non-reliable eigenvalues) will go to zero when $N$ goes to infinity!

\bigskip

{\bf Acknowledgment}. The author would like to thank Professor Huiyuan Li for producing
 the two graphs in the paper.

\end{document}